\documentclass [11pt]{article}

\title {A Faithful Representation of the Singular Braid Monoid\\ 
on Three Strands}
\author {Oliver T. Dasbach
       \thanks{supported by the Deutsche Forschungsgemeinschaft (DFG)}
        \thanks{e-mail: \sl kasten@math.columbia.edu,
       \sl http://www.math.uni-duesseldorf.de/$\sim \!\!$ kasten}\\
          {\small \sl Columbia University}\\
          {\small \sl Department of Mathematics}\\
          {\small \sl New York, NY 10027}\\\\
         Bernd Gemein
        \thanks{e-mail: \sl gemein@rz.uni-duesseldorf.de}\\
        {\small \sl Heinrich-Heine-Universit\"at}\\
          {\small \sl Mathematisches Institut}\\
          {\small \sl Universit\"atsstr. 1, D-40225 D\"usseldorf}
          }

\date{}

\input {epsf}
\usepackage{amssymb}
\usepackage{latexsym}

\newcommand {\se} {{\sigma _1}}
\newcommand {\sz} {{\sigma _2}}
\newcommand {\s} {{\sigma}}
\renewcommand {\t} {{\tau}}
\renewcommand {\O} {{\mathcal O}}
\renewcommand {\aa} {{a_{21}}}
\newcommand {\ab} {{a_{32}}}
\newcommand {\ac} {{a_{31}}}
\newcommand {\mat}[4] 
 {\left (\! \begin{array}{cc}
      #1 & #2\\
      #3 & #4\\
   \end{array} \!
   \right) }
\newcommand {\pr}[2]
   {\langle #1 \vert \, #2 \rangle }

\newcommand {\Z}{\mathbb{Z}}
\newcommand {\Q}{\mathbb{Q}}
\newtheorem{theorem}{Theorem}[section]
\newtheorem{lemma}[theorem]{Lemma}

\newtheorem{remark}[theorem]{Remark}
\newtheorem{proposition}[theorem]{Proposition}
\newtheorem{corollary}[theorem]{Corollary}

\newenvironment {proof}{\noindent {\bf Proof }}{$\Box$\bigskip \medskip}
 
\begin{document}

% \subjclass{57M25, 20F10, 20F36}

\maketitle

\begin{abstract}
We show that a certain linear representation of the singular braid monoid $SB_3$
is faithful. Furthermore we will give a second - group theoretically motivated -
solution to the word problem in $SB_3$.
\end{abstract}

\section {Preliminaries}

While knots are important for our daily life - ties, shoelaces etc. - singular knots
are not. However, since the theory of Vassiliev invariants started (see e.g. \cite{BL}) 
mathematicians became more and more interested in singular knotted objects, i.e.
objects having a finite number of transversal self intersections.
One of these objects are singular braids which form 
the singular braid monoid $SB_n$, i.e. the monoid
generated by the standard generators $\s_1, \dots, \s_{n-1}$ 
of the braid group $B_n$ plus the additional singular generators $\t_1, \dots,
\t_{n-1}$.

A presentation for the singular braid monoid in terms of these generators, 
build up from the usual choice 
of a presentation for the braid group is not hard to deduce 
(cf. \cite{Birman2}).
 
The theory of Vassiliev invariants suggests a homomorphism from $SB_n$
into the integral group ring $\Z B_n$ of the braid
group. This homomorphism is given by mapping a singular generator
$\t_i$ to $\s_i - \s_i^{-1}$ and $\s_i$ to itself (cf. \cite{Birman2}).                 
A famous conjecture of Joan Birman (\cite {Birman2}, \cite{FRZ}) asserts that this 
homomorphism is injective.

As it appeared to be very useful for knot theorist to have a complete understanding of
the braid group $B_3$ (see e.g. \cite{BM}) the purpose of this paper is to give
some analogous results for the braid monoid $SB_3$.

By results of Fenn et al. \cite{FKR} the monoid $SB_n$ embeds in a group $SG_n$.
Using this result we will prove that an extension to $SB_n$ of the Burau representation
for $B_n$ that was defined in \cite{Gemein1} is faithful for $n=3$.
Since the Burau representation of $B_n$ is known to be unfaithful for
$n \geq 6$ \cite{LP} this result cannot carry over at least for $n \geq 6$.
So we will give a second solution to the word problem for $SB_3$ by using the
group theoretical structure of $SG_3$.

It is worth mentioning that the faithfulness of the singular Burau representation
yields an algorithm for the solution of the word problem for $SB_3$ that has a time
complexity which is linear in the length of  the word.   

This paper was written during a visit of the first author
at the Columbia University. He would like to thank Columbia 
for the warm hospitality and especially Joan Birman for many discussions and 
fruitful comments on an earlier version of this text. 

The second author would like to thank Wilhelm Singhof for many useful suggestions and remarks. 

Furthermore both authors are very grateful to E. Mail for her invaluable help.

\section{A new presentation for the singular braid monoid}
We recall the well-known presentation for the monoid of singular braids
on $n$ strands:

\begin{proposition}[Baez \cite {Baez}, Birman \cite{Birman2}]
The monoid $SB_n$ is generated by the elements $$\s_i^{\pm 1}, \, i=1, \dots, n-1, \,
\t_i, \, i=1, \dots, n-1$$ (see Figure \ref{generators}) satisfying the following
relations:
\begin{eqnarray}
\s_i \s_i^{-1}&= & 1 \mbox{ for all } i\\
\s_i \s_{i+1} \s_i = \s_{i+1} \s_i \s_{i+1}
&\mbox{and}& \s_i \s_j = \s_j \s_i \quad \mbox{ for } \, \, j>i+ 1\\
\t_{i+1} \s_i \s_{i+1} &=& \s_i \s_{i+1} \t_i\\
\s_{i+1} \s_i \t_{i+1} &=& \t_i \s_{i+1} \s_i\\
\s_i \t_i &=& \t_i \s_i\\
\s_i \t_j &=& \t_j \s_i \quad \mbox{ for } \,\, j>i+ 1\\
\s_j \t_i &=& \t_i \s_j \quad \mbox{ for } \,\, j>i+ 1\\
\t_i \t_j &=& \t_j \t_i \quad \mbox{ for } \,\, j>i+ 1.
\end{eqnarray}
\end{proposition}

We use this proposition to derive a new presentation for the
singular braid monoid which is more suitable for our purposes:

\begin{figure}[h] 
\begin{center}
\parbox{3.5cm}{\epsfbox{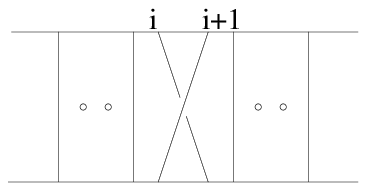}} 
\parbox{3.5cm}{\epsfbox{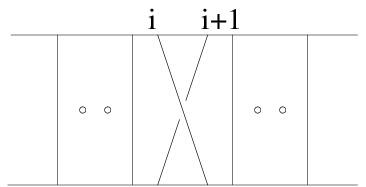}}
\parbox{3.5cm}{\epsfbox{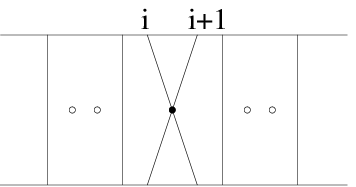}}
\end{center}
\caption{\label{generators} $\sigma_i, \sigma_i^{-1}$ and $\tau_i$}
\end{figure}

\begin{proposition}
The monoid $SB_n$ is generated by the elements $$\s_i^{\pm 1}, \, i=1, \dots, n-1, \, \mbox{ and }
\t$$ satisfying the following relations:
\begin{eqnarray}
\s_i \s_i^{-1}&= &1 \mbox{ for all } \, i  \label{Braid relation 1}\\
\s_i \s_{i+1} \s_i = \s_{i+1} \s_i \s_{i+1}
&\mbox{and}& \s_i \s_j = \s_j \s_i \quad \mbox{ for } \, \, j>i+ 1
\label{Braid relation 3}\\
(\s_2 \s_1)^3 \t &=& \t (\s_2 \s_1)^3\\
\s_i \t &=& \t \s_i \quad \mbox{ for } \, i \neq 2 \label{commutator}\\
\s_2 \s_3 \s_1 \s_2 \t \s_2 \s_3 \s_1 \s_2 \t &=&
\t \s_2 \s_3 \s_1 \s_2 \t \s_2 \s_3 \s_1 \s_2   \mbox{ for } n >  3. \label{t1t3}
\end{eqnarray}
\end{proposition}

\begin{proof}
We will see that the new relations are included in the old ones.
So we will concentrate ourselves to show that all other relations can be derived from
our relations.

First we note that we have
\begin{eqnarray} \label{t_i Definition}
\t_{i+1} = \s_i \s_{i+1} \t_i \s_{i+1}^{-1} \s_i ^{-1}.
\end{eqnarray}

\begin{enumerate}
\item {\bf Relations of the form $\s_i \t_i = \t_i \s_i$}

We will show that these relations can be deduced from $\s_1 \t_1 = \t_1 \s_1$,
Relation (\ref{t_i Definition})
and the relations in the braid group (\ref {Braid relation 1})-
(\ref {Braid relation 3}). We will proceed by induction on $i$:
\begin{eqnarray*}
\s_i \t_i &=& \s_i \s_{i-1} \s_i \t_{i-1} \s_i^{-1} \s_{i-1}^{-1}\\
&=& \s_{i-1} \s_i \s_{i-1} \t_{i-1} \s_i^{-1} \s_{i-1}^{-1}\\
&=& \s_{i-1} \s_i \t_{i-1} \s_{i-1} \s_i^{-1} \s_{i-1}^{-1}\\
&=& \s_{i-1} \s_i \t_{i-1} \s_i^{-1} \s_{i-1}^{-1} \s_i\\
&=& \t_i \s_i
\end{eqnarray*}

\item {\bf Relations of the form $\s_i \t_j = \t_j \s_i, \, \vert i-j \vert >1$}

We will show that all relations of this form can be derived from
the relations $\s_i \t_1 = \t_1 \s_i, \quad i\neq 2$
as well as (\ref{t_i Definition}) and the braid group relations
(\ref {Braid relation 1}) - (\ref {Braid relation 3}).

{\bf Case 1: }
If $1 < j <i-1$ then we
have
\begin{eqnarray*}
\s_i \t_j &=& \s_i \s_{j-1} \s_{j} \t_{j-1} \s_j^{-1} \s_{j-1}^{-1}\\
&=& \s_{j-1} \s_{j} \t_{j-1} \s_j^{-1} \s_{j-1}^{-1} \s_i\\
&=& \t_j \s_i
\end{eqnarray*}
where we use (\ref{t_i Definition}), (\ref {Braid relation 3})
and induction on $j$.

{\bf Case 2: }
In the braid group with $w:=\s_{i+1} \s_{i+2} \s_i \s_{i+1}$ it holds:
\begin{eqnarray}
\s_i w&=& w \s_{i+2} \label{computation 1}.
\end{eqnarray}

Thus, for $j = i+2$:
\begin{eqnarray*}
\s_i \t_{i+2} &=& \s_i w \t_i w^{-1}\\
&=& w \s_{i+2} \t_i w^{-1}\\
&=& w \t_i \s_{i+2} w^{-1}\\
&=& w \t_i w^{-1} \s_i\\
&=& \t_{i+2} \s_i
\end{eqnarray*}
where we use (\ref{t_i Definition}), (\ref{computation 1}) and Case $1$.

{\bf Case 3: } If $j \geq i +3$ then we can proceed as in Case 1 to show that
$$\s_i \t_j = \t_j \s_i, \quad i>1+j$$ follows from our relations and Case $2$.

\item {\bf Relations of the form $\t_i \t_j = \t_j \t_i, \, i<j, j-i>1$}

Our aim is to reduce these relations to $\t_1 \t_3 = \t_3 \t_1$
with the help of (\ref{t_i Definition})
and the relation $\s_i \t_j= \t_j \s_i, \,$ for $\vert j-i \vert \neq 1$.
We use induction on the pair $(i,j)$ that is ordered lexicographically.

If $i>1$ then
\begin{eqnarray*}
\t_i \t_j &=& \s_{i-1} \s_{i} \t_{i-1} \s_{i}^{-1} \s_{i-1}^{-1} \t_j\\
&=& \t_j \s_{i-1} \s_{i} \t_{i-1} \s_{i}^{-1} \s_{i-1}^{-1}\\
&=& \t_j \t_i.
\end{eqnarray*}

In the same way one can deduce $\t_1 \t_j = \t_j \t_1$ for $j>3$.

\item {\bf Relations of the form
$\s_{i+1} \s_i \t_{i+1} = \t_i \s_{i+1} \s_i$}

We will show that these relations follow from Relation (\ref {t_i Definition})
together with the relations in the braid group and the
initial relation:
\begin{eqnarray}
\s_2 \s_1 \t_2 = \t_1 \s_2 \s_1 \label{Induction relation}
\end{eqnarray}
 
We will need the 
following two relations in the braid group which
can be easily tested:
\begin{eqnarray}
\s_i^{-1} \s_{i+1}^{-1} \s_{i-1}^{-1} \s_i^{-1}&=&
\s_{i-1}^{-1} \s_i^{-1} \s_{i+1}^{-1} \s_{i-1}^{-1} \s_i^{-1} \s_{i+1}\\
\mbox{and} \qquad \s_{i-1}^{-1} \s_i^{-1} \s_{i+1}^{-1} \s_{i-1} \s_i \s_{i+1}&=&
\s_i \s_{i+1} \s_{i-1}^{-1} \s_i^{-1}
\end{eqnarray}

\begin{eqnarray*}
\s_i^{-1} \s_{i+1}^{-1} \t_i \s_{i+1} \s_i &=&
\s_i^{-1} \s_{i+1}^{-1} \s_{i-1}^{-1} \s_i^{-1} \t_{i-1}
\s_i \s_{i-1} \s_{i+1} \s_i\\
&=& \s_{i-1}^{-1} \s_i^{-1} \s_{i+1}^{-1} \s_{i-1}^{-1}
\s_i^{-1} \s_{i+1} \t_{i-1} \s_{i+1}^{-1} \s_i \s_{i-1}
\s_{i+1} \s_i \s_{i-1} \\
&=& \s_{i-1}^{-1} \s_i^{-1} \s_{i+1}^{-1} \s_{i-1}^{-1}
\s_i^{-1} \t_{i-1} \s_i \s_{i-1}
\s_{i+1} \s_i \s_{i-1}\\
&=& \s_{i-1}^{-1} \s_i^{-1} \s_{i+1}^{-1} \s_{i-1}
\s_i \t_{i-1} \s_i^{-1} \s_{i-1}^{-1}
\s_{i+1} \s_i \s_{i-1}\\
&=& \s_{i-1}^{-1} \s_i^{-1} \s_{i+1}^{-1} \s_{i-1}
\s_i \s_{i+1} \t_{i-1} \s_{i+1}^{-1} \s_i^{-1} \s_{i-1}^{-1}
\s_{i+1} \s_i \s_{i-1}\\
&=&\s_i \s_{i+1} \s_{i-1}^{-1} \s_i^{-1} \t_{i-1}
\s_i \s_{i-1} \s_{i+1}^{-1} \s_i^{-1}\\
&=&\t_{i+1}
\end{eqnarray*}
\end{enumerate}
Now we are left with the braid group relations, Relation (\ref{commutator})
as well as Relation (\ref {t_i Definition}) and
\begin{eqnarray}
\t_2 &=& \s_1^{-1} \s_2^{-1} \t_1 \s_2 \s_1 \label{left1}\\
\t_1 \t_3 &=& \t_3 \t_1 \label{left2}.
\end{eqnarray}

Relation (\ref{left2}) is equivalent to Relation (\ref{t1t3}).

Relation (\ref{left1}) is equivalent to:
\begin{eqnarray*}
\s_2 \s_1 \s_1 \s_2 \t_1 &=& \t_1 \s_2 \s_1 \s_1 \s_2 \\
\Longleftrightarrow \s_2 \s_1 \s_1 \s_2 \s_1^{2} \t_1 \s_1^{-2} &=&
\t_1 \s_2 \s_1 \s_1 \s_2 \\
\Longleftrightarrow (\s_2 \s_1)^3 \t_1 &=& \t_1 (\s_2 \s_1)^3.
\end{eqnarray*}

Finally, we can skip the Relations (\ref {t_i Definition}) because
the $\t_j, \, j>1$, only occur in these relations and we set $\tau:=\tau_1$.
\end {proof}

\begin{remark} As the braid group $B_n$ admits a presentation with two generators
$\s_1$ and $A:=\s_1 \cdots \s_{n-1}$
for every $n$ (cf. Artin's initial paper \cite{Artin}) 
we can rewrite the presentation for the singular braid monoid in terms of
three generators $\s_1, A$ and $\t_1$. We will omit the 
details here.
\end{remark}

\begin{corollary} \label{Corollary_SB_3}
The monoid $SB_3$ is generated by the elements $\se^{\pm 1}, \sz^{\pm 1}, \t_1$
satisfying the following relations:
\begin{enumerate}
\item $\se \se^{-1} = \sz \sz^{-1}= 1$.
\item $\se \sz \se = \sz \se \sz$
\item $\t_1 (\s_2 \s_1)^3  = (\s_2 \s_1)^3 \t_1$
\item $\se \t_1 = \t_1 \se$
\end{enumerate}
\end{corollary}

The following theorem of Fenn, Keyman and Rourke \cite {FKR} will make our arguments on
$SB_n$ much easier:

\begin{theorem} \label{embedding} 
Let $SG_n$ be the group given by the monoid presentation of $SB_n$ considered as a group
presentation.
Then the natural homomorphism of $SB_n$ into $SG_n$ is an embedding.
\end{theorem}

\section{A faithful representation of $SB_3$}

One can find a representation of the singular braid groups
which is an extension of the famous Burau representation
of the braid groups itself (cf. \cite{Gemein1}). 
For $SB_3$ this representation looks like:

\begin{proposition} 
The map $\beta_s$ given by
$$ \se \mapsto \mat {-t} 1 0 1, \, \sz \mapsto \mat 1 0 t {-t}, \,
\t_1 \mapsto \mat {1 - y -t y} y 0 1$$

yields a representation of the singular braid monoid into a 
matrix ring: $$\beta_s: B_3 \rightarrow M_2(\Z [ t, t^{-1}, y]).$$
\end{proposition}

We will show that this representation is faithful. For this purpose we need the 
following easy consequence (cf. \cite {Fine}) of a theorem of P.M. Cohn 
\cite{Cohn}:

\begin{theorem} \label{Cohn}
Let $d \neq 1,2,3,7,11$ be square-free, i.e. $d$ is not divisible by the
square of an integer,  and $\omega:=\sqrt{-d}$ if $d \equiv 1$ or $2$ modulo $4$ or
$\omega:= (1+\sqrt{-d})/2$ for $d \equiv 3$ modulo $4$
and let $\O_d= \Z + \omega \Z$
be the imaginary quadratic integers in $\Q[\omega]$.
Furthermore, let $A, B$ and $C$ be the following elements in
$PSL_2(\O_d)$:
$$A:= \mat 1 1 0 1, \, B:= \mat 0 {-1} 1 0 \, \mbox{ and } \,
C:= \mat 1 {\omega} 0 1.$$

The subgroup $PE_2(\O_d)$ of $PSL_2(\O_d)$
generated by all matrices of the forms
$$ \mat 1 x 0 1 \,  \mbox{ and } \, \mat 1 0 y 1 $$
for $x$ and $y$ in $\O_d$, has the presentation:

$$PE_2(\O_d) = \pr {A,B,C} {B^2=(A B)^3=[A,C]=1}.$$

With $\Sigma_1:=A, T_1:=C$ and  
$$\Sigma_2:=(A B A)^{-1}=\mat 1 0 {-1} 1$$ we get
- after some easy transformations - the
presentation:
$$PE_2(\O_d) =\pr {\Sigma_1, \Sigma_2, T_1} {(\Sigma_1 \Sigma_2 \Sigma_1)^2=
(\Sigma_1 \Sigma_2)^3=1, \, \Sigma_1 T_1 = T_1 \Sigma_1}.$$
\end{theorem}

Now we have all the tools to show that 
\begin {theorem}
The singular Burau representation $\beta_s: SB_3 \longrightarrow M_2(\Z[t,t^{-1},y])$
is faithful.
\end{theorem}
\begin{proof}
The arguments are essentially the same as in the proof of Magnus and Peluso 
(cf. \cite{Birman3}) 
of the faithfulness of the Burau representation for 3-string braids.
While in this proof the well-known presentation of the group $PSL_2(\Z)$ is used  
we will make use of the presentation of $PE_2(\Z[\omega])$ of Theorem \ref{Cohn} 
for a suitable ring of integers
$\Z[\omega]$. 

For reasons of convenience we choose $d=5$, so $\omega:=\sqrt{-5}$.

We will denote the extension of $\beta_s$ to $SG_3$ also by $\beta_s$.
So $$\beta_s: SG_3 \rightarrow M_2(\Z[t,t^{-1}, y, 1/(1-y-t y)]).$$
The image of $\beta_s$ in this matrix group is - of course - a quotient of $SG_3$.
Furthermore by setting $t:=-1$ and $y:=\omega$ it naturally maps onto $PE_2(\Z[\omega])$.
 
Thus we have a homomorphism from $SG_3$ onto $PE_2(\Z[\omega])$ given by
$\s_1 \mapsto \Sigma_1, \s_2 \mapsto \Sigma_2$ and $\tau_1 \mapsto T_1$ where
$\Sigma_1, \Sigma_2$ and $T_1$ are as in Theorem \ref{Cohn}. 
Moreover by comparing the presentations of $SG_3$ and $PE_2(\Z[\omega])$
we see that the kernel of this homomorphism is the normal closure of 
the element $(\s_1 \s_2 \s_1)^2$, which is - as is easy to see - 
a central element in $SG_3$.
Hence, the kernel is cyclic. 

Thus, the image of $SG_3$ under $\beta_s$ is isomorphic to $SG_3$ modulo a power of 
$\beta_s((\s_1 \s_2 \s_1)^2)$. Since $(\s_1 \s_2 \s_1)^2$ is mapped by $\beta_s$
to the element $$\mat {t^3} 0 0 {t^3}$$ of infinite order, we see that
$SG_3$ must be isomorphic to its image under $\beta_s$.  

Because of the Embedding Theorem \ref{embedding} the theorem follows.
\end{proof}

\subsection{The Birman-Conjecture for $SB_3$} 

Let $\eta: SB_n \longrightarrow \Z B_n$ be the Birman homomorphism which maps
$\t_1 \mapsto \s_1 - \s_1^{-1}$ and $\s_j$ to itself. 
In \cite{Gemein2} the follwing result is proved:

\begin{theorem} \label{Gemein-Injectivity}
Let $b$ and $b'$ be two braids in $SB_n$ with $\eta(b) = \eta(b')$. 
Then we have $\beta_s(b) = \beta_s(b').$
\end{theorem}

Because the singular Burau representation is faithful for $n=3$ this means
that the Birman homomorphism is also faithful, i.e. the Birman conjecture is
valid for $n=3$. We learned that this result was already obtained by
Antal J\'arai \cite {Jarai}.

It might be interesting, however, for the reader to see how the linearity
of $SB_3$ and the injectivity of $\eta$ relate. Therefore we will  
give a sketchy and informal proof of Theorem \ref{Gemein-Injectivity}. 
For full details the reader is referred to \cite{Gemein2}.

For notational reasons we will consider the special case $SB_3$. However, 
the argument holds also for higher $n$. Since the 
determinant of the Burau matrix of a given braid with $m$ singularities 
equals $t^k \cdot (1-y-ty)^m$ for some $k \in \Z$, two 
braids with a different number of singularities cannot map to the same 
matrix under $\beta_s$. The same holds for the Birman homomorphism.
Therefore we may restrict ourselves to braids with a fixed number $m$ of 
singularities. The set of all braids with exactly $m$ singularities will 
be denoted by $SB_3^{(m)}$ in the sequel. 
\par The main idea of the proof is to imitate the Birman homomorphism 
on the level of Burau matrices. If we deal with braids having exactly one 
singularity, this can be done easily: Substituting $y$ by 1 corresponds to a 
right-handed resolution of the singularity, substituting $y$ by $x^{-1}$ 
corresponds to a left-handed resolution of the singularity. This induces a well 
defined homomorphism from the matrix ring $M_2(\Z[t,t^{-1},y])$ into the group 
ring $\Z[M_2(\Z[t,t^{-1}])]$ imitating the Birman homomorphism.
\par If the number of singularities is greater than one, we cannot proceed in the same 
way. The two indicated substitutions would correspond to a right-handed (resp. left-handed) 
resolution of all the singularities. Unfortunately, we also have to consider cases 
where some singularities are resolved in a right-handed way while others are resolved 
in a left-handed way. 
\par Therefore we change our point of view slightly: First, we number the 
singularities of our braid $b$ from 1 to $m$. Afterwards we assign to the $j$-th 
singularity of the braid the matrix
\[ 
\mat {1-y_j-ty_j}  {y_j} 
                0 1 
\]  
rather than the usual matrix  
\[ \mat {1-y-ty} y 0 1 . \]

In this way we define a modification of the Burau matrix of $b$. Its entries take 
values in the polynomial ring $\Z[t,t^{-1},y_1, \dots , y_m]$. Of course, the numbering 
of the singularities was somewhat arbitrary. Therefore we shall regard this modified 
matrices only up to permutation of the indices of the $y_i$. (To be more precise, we 
let the symmetric group $\Sigma_m$ act on $M_2(\Z[t,t^{-1},y_1, \dots , y_m])$ and 
consider the orbits. By abuse of notation we shall denote the set of these orbits also 
with $M_2(\Z[t,t^{-1},y_1, \dots , y_m])$.) It is obvious that we obtain the Burau 
representation of $b$ out of its modified Burau matrix by the projection $p$ which sends 
all the $y_i$ to $y$.
\par We have introduced the modified Burau matrix in order to compute the (regular) 
Burau matrices of all possible resolutions. In fact, define a matrix resolution 
$$r: M_2(\Z[t,t^{-1},y_1, \dots , y_m]) \rightarrow M_2(\Z[t,t^{-1}])$$ as a projection where, 
in addition, any $y_i$ is mapped either to $1$ or to $t^{-1}$. The index $\mu(r)$ is 
defined to be the number of $y_i$ which are sent to $t^{-1}$. Clearly, a given resolution 
is not well defined on our orbits. However, taking formal sums over all possible resolutions 
gives a well defined map $\rho$. So, if $M$ is a modified Burau matrix, then 
$$\rho(M) = \bigoplus_{r} (-1)^{\mu(r)} \cdot r(M).$$ Note that the sum in the formula is 
a formal sum in the group ring $\Z[M_2(\Z[t,t^{-1}])].$
\par Easy calculations show that the application $\rho$ corresponds to the Birman 
homomorphism $\eta$ on the level of matrices. In fact, we get the following commutative 
diagram:

\unitlength1cm
\begin{picture}(13,6)
\put(3.2,0){$\Z[B_3]$}
\put(3.3,3){$SB_3^{(m)}$}
\put(8.5,0){$\Z[M_2(t,t^{-1})]$}
\put(8,3){$M_2(t,t^{-1},y_1,\dots,y_m)$}
\put(8.5,5){$M_2(t,t^{-1},y)$}
\put(3.6,2.5){\vector(0,-1){2}}
\put(9.3,2.5){\vector(0,-1){2}}
\put(4.6,0.1){\vector(1,0){3.5}}
\put(4.5,3.1){\vector(1,0){3.2}}
\put(9.3,3.6){\vector(0,1){1}}
\put(4.3,3.6){\vector(3,1){4}}
\put(3.9,1.5){$\eta$}
\put(9.6,1.5){$\rho$}
\put(9.6,4){$p$}
\put(6,0.3){$\Z[\beta]$}
\put(6,2.6){$\tilde \beta$}
\put(5.9,4.5){$\beta_s$}
\end{picture}

\par \bigskip 

\bigskip

Here $\tilde \beta$ denotes the application which maps a braid to its corresponding modified 
Burau matrix. With $\Z[\beta]$ we denote the extension of the usual (regular) Burau 
homomorphism to the group rings.
\par We now claim:
\begin{description}
\item[Claim 1]  Let $M,M' \in M_2(x^{\pm 1}, y_1, \dots , y_m)$ be two elements in the 
image of $\tilde \beta$ with $\rho(M) = \rho (M')$. Then we have $p(M) = p(M')$.
\end{description}
Let us assume that the claim is true. then the proof of Theorem 3.4 becomes easy diagram 
chasing:
\par Let $b,b'$ be elements of $SB_3^{(m)}$ and suppose that $\eta(b) = \eta(b')$. It follows 
that $(\Z[\beta] \circ \eta)(b) = (\Z[\beta] \circ \eta)(b')$ and by commutativity of the 
diagram that
$(\rho \circ \tilde \beta) (b) = (\rho \circ \tilde \beta) (b')$. Using Claim 1 we get 
$(p \circ \tilde \beta) (b) = (p \circ \tilde \beta) (b')$ and - again by commutativity of 
the diagram - $\beta (b) = \beta (b')$.
\par Thus, we only have to show that Claim 1 holds. This is the most technical part of the 
proof. In fact, we have to figure out in how far the matrices of the formal sum $\rho(M)$ 
determine the matrix $M$. This leads to one system of linear equations for each 
entry of the matrix, which 
may be solved after having 
observed the following two facts:
\begin{enumerate}
\item If $M$ is in the image of $\tilde \beta$, then each $y_i$ cannot appear in the matrix 
with powers greater than $1$.
\item We may use the determinant of the matrices in our formal sum in order to compute the 
index of the resolution which has produced them. This fact is important when solving the 
equations.
\end{enumerate}
With these two observations and some tedious computations, we derive that two matrices $M$ 
and $M'$ are mapped to the same formal sum under $\rho$ if their entries differ by 
permutations of the indices of the $y_i$.  
Hence, they vanish under the projection $p$.

\section{A second solution to the word problem in $SB_3$}

To give a second solution to the word problem for $SB_3$, i.e. the problem
whether two elements in $SB_3$ are equivalent,  
we will need Britton's Lemma that can be applied to the group $SG_3$
\begin{lemma}[Britton \cite{Britton}]
Let $H= \pr {S} {R}$ be a presentation of the group $H$ with a set of generators
$S$ and relations $R$ in these generators. 

Furthermore let $G$ be a HNN-extension of $H$ of the following form:
$$ G = \pr {S,t} {R, t^{-1} X_i t = X_i, i\in I}$$ for some index set $I$,
where $X_i$ are words over $S$.

Let $W$ be a word in the generators of $G$ which involves $t$.

If $W=1$ in $G$ then $W$ contains a subword $t^{-1} C t$ or $t C t^{-1}$
where $C$ is a word in $S$, and $C$, regarded as an element of the group $H$,
belongs to the subgroup $X$ of $H$ generated by the $X_i$.
\end{lemma}

We will rather solve the word problem for $SG_3$ than for $SB_3$. 
By Corollary \ref{Corollary_SB_3} $SG_3$ has a presentation as in Britton's Lemma.  
So to solve the word problem in $SG_3$ 
we have to decide whether a given word in the generators of
$B_3$  is element of the subgroup $H_3$ generated by the elements with which $\t_1$ 
commutes:
$\s_1$ and $(\s_2 \s_1)^3$.
  
This decision problem, called membership problem, would not be hard to solve with
the help of the Burau representation of $B_3$. However we promised to give a
puristic proof which gives more hope for a generalization to braids and singular
braids with more than three strands.

So we will choose the approach of Xu \cite{Xu}   for the word problem for $B_3$ - which
was generalized most recently to arbitrary $B_n$ by Birman, Ko and Lee \cite{BKL} -
to solve the membership problem for the subgroup $H_3$.

We briefly recall this  approach using the notation of Birman, Ko and Lee.
The first step is to rewrite the presentation of $B_3$ in terms of the new
generators: $\aa:=\s_1, \ab:=\s_2$ and $\ac:=\s_2 \s_1 \s_2^{-1}$. 

So we get a new presentation      
$$
B_3 = \pr {\aa, \ab, \ac}{\ab \aa = \ac \ab = \aa \ac}.
$$ 

Using the element $\delta := \ab \aa$ one can show now that every element
of $B_3$ can be brought into a unique normal form $\delta ^k P$ for some $k$ with
$P$ a positive word, i.e. only positive exponents
occurs, in the generators $\aa, \ab$  and $\ac$, such that none of the subwords
$\ab \aa$, $\ac \ab$ or $\aa \ac$ appear in $P$.   

\begin{lemma} \label{membership_problem}
The membership problem for the subgroup $H_3$ of $B_3$ generated by the elements
$\s_1$ and $(\s_2 \s_1)^3$ can be solved.
\end{lemma} 

\begin{proof}
First of all we see that $H_3$ is abelian and $(\s_2 \s_1)^3=\delta^3$. Therefore 
if we want to bring a word into the normal form we only have to look for $\s_1^k$ for 
$k \in \Z$.
If $k$ is not negative, then the normal form for an element $(\s_2 \s_1)^{3l} \s_1^k$
is simply $\delta ^{3l} \aa ^k$.

If $k\leq 0$ then the following identities are easy to see: 
\begin{eqnarray*}
(\s_2 \s_1)^{3l} \s_1^{3k} &=& \delta^{3l+3k} (\ac \aa \ab)^{-k}\\
(\s_2 \s_1)^{3l} \s_1^{3k-1} &=& \delta^{3l+3k-1} \ab (\ac \aa \ab)^{-k}\\
(\s_2 \s_1)^{3l} \s_1^{3k-2} &=& \delta^{3l+3k-2} \aa \ab (\ac \aa \ab)^{-k}.
\end{eqnarray*}

Therefore for every word $w$ in the braid group $B_3$ we can bring it into its
unique normal form and compare this form with the normal forms for the elements
in $H_3$. Hence, the membership problem for $H_3$ is solvable.
\end{proof}

Thus we have proved:
\begin{theorem} 
Given two words $w_1 = b_1 \t_1 b_2 \t_1 \cdots \t_1 b_m$ and 
$w_2= c_1 \t_1 c_2 \t_1 \cdots \t_1 c_l$ 
in the generators $\s_1, \s_2$ and $\t_1$ in $SB_3$, where the $b_j$ and $c_j$ are words in $B_3$.

Then $w_1$ and $w_2$ are equal in $SB_3$ if and only if $b_m c_l^{-1}$ is in the subgroup $H_3$ of $B_3$
that is generated by $\s_1$ and $(\s_2 \s_1)^3$ and 
$b_1 \t_1 b_2 \t_1 \cdots b_{m-1} b_m$ and $c_1 \t_1 c_2 \t_1 \cdots c_{l-1} c_l$ are equal in $SB_3$.

This gives a solution to the word problem in $SB_3$ because the membership problem for $H$ is solvable by
Lemma \ref{membership_problem} and the word problem in $B_3$ is solvable - as mentioned 
above - as well.
\end{theorem}

\subsection{An algebraic proof of the embedding theorem for $SB_3$ into $SG_3$} 

Actually - as a Corollary of our approach - one can get a purely algebraic proof
of the Embedding Theorem \ref{embedding} of \cite{FKR} for the special case $n=3$:

\begin{corollary}
$SB_3$ embeds into $SG_3$.
\end{corollary}

\begin{proof}
We have to show that if two elements $w_1$ and $w_2$ in $SB_3$ are different then
their images in $SG_3$ are different. 
By the HNN-structure of $SG_3$ the subgroup $B_3$ embeds in it. Furthermore 
two elements $w_1$ and $w_2$ that map to the same element in $SG_3$ 
must have the same number of singular points. 

Now let $w_1 = b_1 \t_1 b_2 \t_1 \cdots \t_1 b_m$ and 
$w_2= c_1 \t_1 c_2 \t_1 \cdots \t_1 c_m$, $b_j$ and $c_i \in B_3$,  be two different elements 
of $SB_3$ that map
to the same element in $SG_3$, by slight abuse of notation also denoted by the same
word. 
We assume that $w_1$ and $w_2$ are minimal examples with respect to the number of
singular points.
Since $w_1 \beta \neq w_2 \beta  \iff w_1 \neq w_2$ for an element $\beta \in B_3$
we further may assume that $c_m = 1$. 
Then by Britton's Lemma $b_m$ must lie in the subgroup $H_3$ defined above. Since all
the elements of $H_3$ commute with $\t_1$ we have $\t_1 b_m = b_m \t_1$ both in $SB_3$ 
and $SG_3$.

So $w_1$ is equal within $SB_3$ to $w_1 = b_1 \t_1 b_2 \t_1 \cdots b_{m-1} b_m \t_1$.
Now consider the two word 
$w_1' = b_1 \t_1 b_2 \t_2 \cdots b_{m-1} b_m$ and  
$w_2'= c_1 \t_1 \cdots c_{m-1}$ in $SB_3$.
These two words represent different elements in $SB_3$ - otherwise 
we would have $w_1 = w_1' \t_1 = w_2' \t_1 = w_2$ - but map to the same element in $SG_3$.
This contradicts our assumption. 
\end{proof}

\providecommand{\bysame}{\leavevmode\hbox to3em{\hrulefill}\thinspace}

\end{document}